\documentclass[11pt,letterpaper]{article}
\usepackage[margin=1in]{geometry}
\usepackage[T1]{fontenc}
\usepackage[utf8]{inputenc}
\usepackage{lmodern}
\usepackage{amsmath,amssymb,amsfonts,amsthm,bm}
\usepackage{booktabs,microtype,graphicx,float}
\usepackage[backend=biber,style=authoryear,natbib=true,sorting=nyt,maxcitenames=2,maxbibnames=99]{biblatex}
\usepackage{xcolor}
\usepackage{hyperref}
\hypersetup{colorlinks=true,linkcolor=blue,citecolor=blue,urlcolor=blue,hypertexnames=false}
\usepackage[capitalise,noabbrev]{cleveref}
\numberwithin{equation}{section}
\theoremstyle{plain}
\newtheorem{theorem}{Theorem}
\newtheorem{lemma}{Lemma}
\newtheorem{proposition}{Proposition}
\newtheorem{corollary}{Corollary}
\theoremstyle{definition}

\newtheorem{assumption}{Assumption}

\newcommand{\bbE}{\mathbb{E}}
\newcommand{\bbP}{\mathbb{P}}
\newcommand{\bbR}{\mathbb{R}}
\newcommand{\calN}{\mathcal{N}}
\newcommand{\Sym}{\mathfrak{S}}

\DeclareMathOperator{\Var}{Var}
\DeclareMathOperator{\Cov}{Cov}

\title{\bfseries Noise Calibration and Exact Recovery for\\ Multiview Ordinal Rearrangement Data}
\author{Fanghua Chen \and Yizhou Cai \and Ting Fung Ma \and Yun Yang}
\date{October 2026}
\begin{document}
\maketitle
\begin{abstract}

We studied simultaneous recovery of a latent group count, subject partition, and signature means from multiview ordinal ranks. A normalized reference rearrangement score is an affine Spearman correlation; the problem concerns the information retained by several such scores and their noise geometry. We develop an ordered-gap estimator with direction-specific concentration radii. Its finite-sample guarantee
separates signature noise, missing groups, and candidate coverage, and allows the group count to grow under an attainable-margin condition. A coordinatewise partition followed by label-tuple refinement recovers groups even when every individual coordinate has population collisions. A binary-code construction explains why the one-dimensional inverse-group margin restriction need not govern a retained coordinate bank. We derive actual-law concentration from bounded ranking modules and a fixed-dimensional combinatorial central limit theorem under an observed permutation representation. A binary ordinal submodel gives logarithmic exact-recovery information scales. A subject-effect construction distinguishes persistent overlap from reducible measurement noise, and an additional ordinal anchor
quantifies how new information lowers that error floor. Five focused simulation scenarios illustrate direction-specific calibration, growing-group recovery, coordinate collision resolution, variance reduction, and the effect of an additional ordinal anchor. The simulations are presented with their main specifications and interpretation limits.
\end{abstract}
\noindent\textbf{Keywords:} ordinal ranks; multiview clustering; reference signatures; noise calibration; exact partition recovery; combinatorial central limit theorem.

\section{Introduction}
\label{p1:sec:introduction}

We study the recovery of a subject partition from multiple ordinal profiles, with possibly different numbers of items in different views. The signal is an ordering pattern relative to specified reference profiles. Raw response levels are deliberately discarded. This invariance is appropriate when changes of scale are nuisance effects, but it also changes the inferential target: groups distinguished only
by features removed by ranking cannot be recovered from ranks.

Probabilistic ranking models provide the natural comparison class.
Distance-based models originate with \citet{Mallows1957} and
\citet{FlignerVerducci1986}; \citet{MurphyMartin2003} developed mixture
models for heterogeneous ranking populations. Multivariate and partial
rankings are already covered by model-based methods
\citep{JacquesBiernacki2014,Rankcluster2014}. More recent computational
developments include mixtures with Spearman distance and their
implementation in MSmix \citep{CrispinoEtAl2023,MSmix2025}.
Thus neither ranking mixtures nor multiview ranking clustering is new.
Weighted distance-based ranking mixtures also precede this work
\citep{LeeYu2012,QianYu2019}. Their item or position weights belong to
the ranking likelihood, whereas our simplex weights combine
reference-signature coordinates. Provable learning and identifiability
of Mallows mixtures have also been studied \citep{AwasthiEtAl2014}.

Our construction uses normalized rearrangement scores relative to a
reference in each view. Its connection with order-statistics
association measures \citep{WeichaoXu2007} is explicit. For complete
ranks, the signed score is exactly Spearman's rank correlation, and the
unsigned score is an affine transform of squared Spearman distance.
We therefore do not claim a new univariate correlation coefficient.
The statistical question is what a finite collection of these
reference features preserves, how their noise behaves, and when one
weighted coordinate suffices for an unknown subject partition.

Two logically different issues arise. Distinct signature means admit a
separating simplex projection outside a finite collection of hyperplanes.
This geometric fact does not imply separation of noisy observations,
uniqueness of the weight, or identification of the number of groups.
Moreover, different ranking distributions can have identical mean
signatures for every reference. We give an explicit example and restrict
mean-based recovery claims to mean-identifiable active classes.

The central estimator is an ordered-gap rule with a calibrated noise
radius. Sorting and gap connectivity are classical clustering devices;
our contribution is their actual-law calibration for ordinal reference
features and the resulting joint recovery analysis. This target differs
from density cluster-tree recovery \citep{ChaudhuriEtAl2014ClusterTree}:
the active groups are specified by a latent model and verified noise
scales, rather than by connected components of density level sets. For each candidate weight, it sorts the projected scores and
cuts only gaps larger than the noise threshold. It then chooses the
largest certified group count among the candidate projections.
On a uniform concentration event, no projection can split a true group;
one sufficiently separated projection prevents all merging. This yields
joint recovery of the group count and the partition without assuming
that a multi-penalty population loss already has the desired minimizer.
The weights themselves need not be uniquely identified.

The analysis separates the numbers of subjects, views, and within-view
units, denoted by $m$, $T$, and $n_t$. The permutation length for a
view-specific statistic is $n_t$, whereas $m$ controls class occupancy
and simultaneous recovery. Recovery bounds allow $T$ and the active
group count to vary when their stated noise and gap conditions hold.
The multivariate Gaussian limit is stated for fixed dimension; it is
not used as a high-dimensional likelihood approximation.

Our concentration conditions concern the actual latent ranking law.
An auxiliary uniform-permutation reference supplies a combinatorial
calibration, not a model for the observed alternatives. We prove an
actual-law Bernstein bound for independent bounded ranking modules,
give a complete-rank pair-swap family with an explicit variance of
order $n_t^{-1}$, and state a combinatorial CLT under an observed-law
permutation linearization
\citep{hoeffding1951combinatorial,Tuvaandorj2024StratifiedCLT}.
Shared subject effects can instead leave a positive noise floor; a
second ordinal construction shows why increasing the item counts need
not then yield exact recovery.

The theory provides a quantitative simplex-search guarantee, an oracle
cutoff bound governed by projected variance, and finite-sample joint
recovery for the calibrated estimator. In a binary ordinal submodel,
an information lower bound matches the logarithmic subject-count
dependence of the sufficient recovery scale, up to constants. These
results identify the roles of representation, weight search, and
irreducible noise separately. A bounded one-dimensional score has a
minimum-gap upper bound of $1/(K_0-1)$; a positive separating direction
alone therefore gives no uniform recovery guarantee as the number of
groups grows. We give additional sufficient rates under an explicit
attainable-margin condition, and distinguish a conservative universal
random pool from a bank that actually attains that margin.

A one-dimensional score in $[0,1]$ has minimum gap at most
$1/(K_0-1)$, but this is a restriction of compression, not a universal
information bound for the multiview problem. We therefore also analyze
a feasible coordinatewise gap partition and its refinement by tuples
of coordinate labels. This estimator uses neither population centers
nor the active count. A resolved-level condition yields exact recovery
despite collisions in every individual coordinate. In a binary-code
model, the required measurement budget depends on the number of code
coordinates and their replication, rather than on an unavoidable
$K_0^2$ factor. This conclusion is specific to the stated code geometry.

The computational task of learning candidate directions is kept separate
from the recovery theorem. Neural and bootstrap-based generators can
provide candidates, but their usefulness must be established through a
coverage probability. The present paper concentrates on representation,
concentration, and exact recovery for an explicit finite pool. Learning
algorithms and their computational validation form a separate project.

Section \ref{p1:sec:method} defines the ordinal representation and the
calibrated estimator. Section \ref{p1:sec:theory} gives the assumptions,
recovery results, and limitations. Section \ref{p1:sec:simulation}
reports five focused scenarios on noise calibration, growing group
counts, coordinate collisions, variance geometry, and additional ordinal
information. No fitted ranking-mixture comparison or neural-training
claim is made. Technical proof details are not included in this preprint version.

%------------------------------------------------------------
% Methodology
%------------------------------------------------------------

\section{Methodology}
\label{p1:sec:method}

\subsection{Ordinal observations and the recovery target}
\label{p1:subsec:multitensor}

Let $\mathcal D=\{\mathcal X_1,\ldots,\mathcal X_m\}$, where
$\mathcal X_r=(X_r^1,\ldots,X_r^T)$ and $X_r^t\in\Sym_{n_t}$ is the
vector of ordinal ranks attached to the $n_t$ item identities in view
$t$. A permutation is identified with its rank vector. If numerical
responses are available, ranking is performed within each view and the
resulting ranks are returned to their original item positions. Sorting
both profiles and then discarding the item identities would remove the
alignment signal.

Throughout, $m\geq1$ and $n_t\geq2$.
The primary theory uses complete ranks. Ties require a stated
preprocessing rule. Independent random ordering within tied blocks
produces complete ranks, with its extra randomness included in the
observed law; it does not preserve all information in the original ties.
Midranks define a different branch whose denominator and concentration
must be checked separately. A constant profile is flagged as
uninformative before preprocessing. We do not define an undefined
normalized score by assigning it the value $1/2$.

Write $Y_r\in\{1,\ldots,K_0\}$ for the active latent group and
$\pi_k=\bbP(Y_r=k)$. The observations are independent across subjects;
the entries within a view and the views of one subject may be dependent.
The joint ordinal law given $Y_r=k$ is $P_k$, and its view marginals
need not determine its dependence structure. A group may represent a
joint pattern of view-specific types. The target is the partition
induced by these ordinal laws, subject to the signature separation
condition below, rather than recovery of arbitrary raw-data mixtures.

The notation $m$ always counts subjects, $T$ counts views, and $K_0$
counts active groups. Put $n_{\min}=\min_t n_t$. Their asymptotic
relationship is specified by the recovery conditions rather than by
a fixed aspect ratio.

\subsection{Reference signatures and information retained}
\label{p1:subsec:tensor_signature}

Let $a^t\in\Sym_{n_t}$ be a fixed reference rank vector. References
estimated on independent calibration data can be treated conditionally;
the stated bounds must then hold conditionally, and calibration failure
probabilities must be added. For the observed item allocation, put
\begin{equation}
Q_r^t=\sum_{j=1}^{n_t}a_j^tX_{rj}^t,\qquad
D_t=\frac{n_t(n_t^2-1)}6.
\label{p1:eq:raw_rearrangement_statistic}
\end{equation}
The two extremal rearrangements are
$$
Q_{\max,t}=\frac{n_t(n_t+1)(2n_t+1)}6,\qquad
Q_{\min,t}=\frac{n_t(n_t+1)(n_t+2)}6.
$$
The signature coordinate and full signature are
\begin{equation}
R_r^t=\frac{Q_r^t-Q_{\min,t}}{D_t},\qquad
R_r=(R_r^1,\ldots,R_r^T)^\top\in[0,1]^T.
\label{p1:eq:tensor_signature}
\end{equation}
The deterministic denominator in \eqref{p1:eq:tensor_signature} removes
the need for a denominator delta method in the complete-rank branch.
The observed allocation is retained in
\eqref{p1:eq:raw_rearrangement_statistic}.

For $L$ references, stack the $LT$ coordinates and use the same
construction with dimension $p=LT$ in place of $T$. Multiple
references are measurements of the same subject, not independent
replicates. Their mean signatures are linear functions of the mean
rank vectors. They can resolve a collision only when an added
reference contributes a population difference.

Let $\eta_k=\bbE(R_r\mid Y_r=k)$. The population mean geometry is
$$
\Delta_\infty=\min_{k<\ell}\|\eta_k-\eta_\ell\|_\infty,\qquad
\gamma(w)=\min_{k<\ell}|w^\top(\eta_k-\eta_\ell)|.
$$
The minimum over an empty pair set is $+\infty$. Recovery based on
these means requires their pairwise distinction; it is not implied by
$P_k\neq P_\ell$ alone.

\subsection{Weights, population cutoffs, and variance geometry}
\label{p1:subsec:wrc_rule}

Let $\Delta_T^{\mathrm w}=\{w\geq0:\mathbf1^\top w=1\}$. Define
$$
Z_r(w)=w^\top R_r,\qquad
G_{w,c}(R_r)=1+\sum_{j=1}^{K-1}\mathbf1\{Z_r(w)>c_j\}.
$$
URC uses $w^{\mathrm U}=T^{-1}\mathbf1$. WRC searches over the simplex.
For two known population means, a weight with
$w^\top(\eta_1-\eta_2)\neq0$ admits the midpoint cutoff
$c(w)=w^\top(\eta_1+\eta_2)/2$.
For $K_0$ distinct projected centers
$\zeta_{(1)}(w)<\cdots<\zeta_{(K_0)}(w)$, the population cutoffs are
$c_j(w)=\{\zeta_{(j)}(w)+\zeta_{(j+1)}(w)\}/2$.

At most $\binom{K_0}{2}$ pair-specific directions suffice as a
pairwise measurement bank when each pair is separable. This is an
upper bound, not a necessary number of directions. In fact, a single
simplex direction generically separates every finite collection of
distinct means. Theorem \ref{p1:thm:separating_weight_exists} quantifies
the search and its dimension dependence.

For actual-law Bernstein parameters $A_k,b_k$, define
\begin{equation}
\begin{gathered}
v(w)=\max_k\|A_k^\top w\|_2^2,\qquad
\beta(w)=\max_k b_k\|w\|_\infty,\\
\mathcal I(w)=\frac{\gamma(w)^2}{8v(w)+4\beta(w)\gamma(w)}.
\end{gathered}
\label{p1:eq:projected_information}
\end{equation}
When the numerator is positive and the denominator is zero, the ratio
is $+\infty$. A zero gap is assigned information zero.
Equation \eqref{p1:eq:projected_information} exposes the noise cost of a
direction. A maximum raw gap need not maximize recovery information.
Pairwise variance-aware directions can be used as pilot proposals,
but latent labels are not supplied to the implementable estimator.

\subsection{The noise-calibrated ordered-gap estimator}
\label{p1:subsec:wrc_unknown_k}

Fix a calibrated radius $\rho\geq0$ and a finite nonempty candidate
collection $\mathcal W=\{w^1,\ldots,w^N\}\subset\Delta_T^{\mathrm w}$.
The radius must control the actual signatures. Calibration from an
auxiliary uniform permutation is valid only if its transfer to that
actual law has been established. Section \ref{p1:subsec:actual_limits}
gives a fully verified ordinal family and a general sufficient
concentration condition.

For each $w$, sort the projected observations as
$z_{(1)}(w)\leq\cdots\leq z_{(m)}(w)$ and define
\begin{equation}
\mathcal B_\rho(w)=
\{i<m:z_{(i+1)}(w)-z_{(i)}(w)>3\rho\},\qquad
\widehat K_\rho(w)=1+|\mathcal B_\rho(w)|.
\label{p1:eq:certified_gap_count}
\end{equation}
Select $\widehat w$ to maximize $\widehat K_\rho(w)$ over
$\mathcal W$. Among maximizers with at least two groups, maximize the
smallest selected gap, denoted by $\widehat g(w)$; remaining ties are
resolved by the fixed ordering of candidates. If all counts equal
one, use the first candidate. Set $\widehat K=\widehat K_\rho(\widehat w)$.
For the selected boundaries $i_1<\cdots<i_{\widehat K-1}$, take
\begin{equation}
\widehat c_j=
\frac{z_{(i_j)}(\widehat w)+z_{(i_j+1)}(\widehat w)}2,\qquad
\widehat Y_r=G_{\widehat w,\widehat c}(R_r).
\label{p1:eq:wrc_cutoff_rule}
\end{equation}
When $\widehat K=1$, the cutoff vector is empty. Thus a single group
is an admissible outcome, without an artificial null component.

This estimator maximizes a certified count, not an uncalibrated
separation loss. On a uniform radius event, a true group has projected
diameter at most $2\rho$ for every simplex weight and therefore
cannot be split by \eqref{p1:eq:certified_gap_count}. A candidate with
$\gamma(w)>5\rho$ separates all groups. These two statements identify
$K_0$ and the partition simultaneously.

The candidate collection may include coordinate directions, URC,
Dirichlet draws, and directions learned from the same data. The
uniform radius event controls even data-dependent directions.
Independent uniform-simplex draws provide a separate probability
guarantee that at least one sufficiently separated candidate is present.
Exact search within a finite pool is inexpensive: constructing and
sorting its scores costs $O(NmT+Nm\log m)$, after the signatures have
been computed. This is a pool guarantee, not a claim of global
optimization over the continuum of weights.

The threshold $3\rho$ is deliberately conservative. It prevents
splitting on the calibrated event but may merge groups that an oracle
with known centers can distinguish. An extension replaces $3\rho$ by
$\vartheta r_j$ for candidate $w^j$, where $\vartheta\geq2$ and $r_j$
controls its projected error simultaneously over subjects. The same
count maximization and gap tie-break define this extended estimator.
The choice $\vartheta=2$ still prevents splitting, because the rule
cuts only gaps strictly larger than $2r_j$. Its sufficient population
gap is $4r_j$, instead of $5r_j$ at $\vartheta=3$. This change improves
the sufficient constant but does not remove the conservatism of a
uniform high-probability noise bound.

For a fixed bank of size $J$, valid deterministic Bernstein envelopes
$\bar v_j\geq v(w^j)$ and $\bar c_j\geq\beta(w^j)$ give radii
$$
r_j(x)=\sqrt{2x\bar v_j}+x\bar c_j,\qquad
x=\log(2mJ/\alpha),\quad 0<\alpha<1.
$$
These radii can exploit cross-view variance reduction, as proved in
Theorem \ref{p1:thm:directional_calibration}. For an independently
pilot-generated bank, the envelopes must hold conditionally on that
pilot; a calibration failure probability is added if necessary.
Estimating a covariance matrix and substituting it without an upper
confidence argument does not certify these radii. An unrestricted
same-data learned bank remains covered by the original uniform
coordinate event, rather than automatically by the finite-bank bound.

\subsection{Combining coordinate partitions}
\label{p1:subsec:tuple_method}

Fix a nonempty bank $w^1,\ldots,w^L$ independently of the recovery data,
and put $Z_r^\ell=(w^\ell)^\top R_r$. Give each coordinate a valid
simultaneous error radius $r_\ell(x)$. In each sorted coordinate,
cut gaps strictly larger than $\vartheta r_\ell(x)$, with
$\vartheta\geq2$, and denote the resulting block label by
$\widehat b_{r\ell}$. Assign the same final label exactly when
two subjects have identical label tuples:
\begin{equation}
\widehat C(r)=\{r':\widehat b_{r'\ell}=\widehat b_{r\ell}
                       \text{ for every }\ell\leq L\},\qquad
\widehat K_{\rm tuple}=|\{\widehat C(r):r\leq m\}|.
\label{p1:eq:tuple_estimator_v2}
\end{equation}
Thus the estimator takes the common refinement of coordinate
partitions. It never supplies $K_0$, a population center, or a true
label to the fitting procedure. Construction and sorting cost
$O(mLT+mL\log m)$ for a dense bank, followed by tuple grouping.
Repeated or constant population coordinates are allowed.
Theorem \ref{p1:thm:tuple_recovery_v2} gives its guarantee.

A second option forms a graph on the subjects, joining $r,r'$ when
$|Z_r^\ell-Z_{r'}^\ell|\leq2r_\ell(x)$ in every coordinate, and
returns its connected components. Proposition
\ref{p1:prop:bank_graph_v2} covers this option under a direct bank-distance
condition. It can use $O(m^2L)$ pair comparisons, whereas the tuple
rule is based on separate sorting. Only the tuple rule is evaluated
in the coordinate-combination experiment below.

\subsection{Candidate generators and likelihood-based comparisons}
\label{p1:subsec:wrc_net}

A neural map or a generalized bootstrap procedure may supply additional
simplex directions. Such generators enter the present theory only through
the probability that their bank contains a sufficiently separated
candidate. No convergence of neural optimization is required for the
independent uniform-pool guarantee. Learning and uncertainty assessment of
these generators are separate computational questions.

For comparison, a conditionally independent multiview Mallows mixture
has joint probability mass
$$
p(x)=\sum_{k=1}^{K}\pi_k
\prod_{t=1}^T
\frac{\exp\{-\theta_{kt}d_t(x^t,\sigma_{kt})\}}
     {\mathcal Z_t(\theta_{kt})}.
$$
This is a parametric model on ranks, unlike the signature conditions.
Our allowance for dependent views is conditional on their verified
concentration, not a guarantee under unrestricted dependence.
For two common-reference Spearman-Mallows components, their posterior
log odds are affine in the reference signatures. Thus the cutoff
family overlaps exactly with a Bayes rule in that specialization.
For general locations and concentrations, one fixed reference
projection need not reproduce the mixture's posterior partition.
Neither universal statistical superiority nor a universal runtime
advantage over EM is claimed.

Gaussian BIC-style scores are possible computational comparators.
They evaluate a working density on ordinal signatures; they are not
validated by a weak CLT. A discrete signature law and a nondegenerate
continuous Gaussian law are mutually singular at finite sample sizes.
Genuine likelihood-order consistency requires a separately specified
likelihood and its own approximation and optimization analysis
\citep{Schwarz1978,Keribin2000}.
The BIC penalty arises from a marginal-likelihood approximation under
model assumptions; it is not a generic unbiased KL-risk penalty.
The present order-recovery theorem instead concerns
\eqref{p1:eq:certified_gap_count}--\eqref{p1:eq:wrc_cutoff_rule}.

%----------------------------------------------------------
%Theory
%----------------------------------------------------------

\section{Theory}
\label{p1:sec:theory}

\subsection{Assumptions}
\label{p1:subsec:theory_assumptions}

The concentration and recovery statements are finite-sample statements.
Their asymptotic consequences allow triangular arrays in $m,T,n_t,K_0$.
Only the multivariate CLT below fixes $T$.

\begin{assumption}[Independent subjects and occupancy]\label{p1:ass:sampling_abundance}
The pairs $(\mathcal X_r,Y_r)$ are iid at each experiment, with
$1\leq K_0<\infty$ and $\pi_{\min}=\min_k\pi_k>0$.
For asymptotic recovery, assume
$$m\pi_{\min}-\log(K_0\vee2)\longrightarrow\infty.$$
\end{assumption}

The occupancy condition is used to recover all active groups, rather
than only those represented in the sample.

\begin{assumption}[Active mean separation]\label{p1:ass:quantitative_identifiability}
For the chosen reference signatures, the active means satisfy
$$\eta_k=\bbE(R_r\mid Y_r=k),\qquad
\Delta_\infty=\min_{k<\ell}\|\eta_k-\eta_\ell\|_\infty>0
\quad(K_0\geq2).$$
\end{assumption}

This restriction is imposed only for positive mean-recovery results.
It is not asserted for every pair of different ordinal distributions.
For multiple references it is imposed on the augmented signature.

\begin{assumption}[Actual-law Bernstein control]\label{p1:ass:signature_concentration}
There exist deterministic $A_k$ and $b_k\geq0$ such that, for every $u$,
$$\log\bbE(e^{u^\top(R_r-\eta_k)}\mid Y_r=k)
\leq\frac{\|A_k^\top u\|_2^2}{2(1-b_k\|u\|_\infty)}
\quad\text{if }b_k\|u\|_\infty<1.$$
Write $a_{kt}=\|A_k^\top e_t\|_2$.
\end{assumption}

The variance proxy is $A_kA_k^\top$. No diagonal covariance assumption
is made. Lemma \ref{p1:lem:bounded_rank_modules} derives this MGF bound
from a primitive independent-module representation.

\begin{assumption}[Optional shrinking coordinate scale]\label{p1:ass:shrinking_noise}
For a constant $C$ uniform over the experiment sequence,
$$a_{kt}\leq Cn_t^{-1/2},\qquad b_k\leq Cn_{\min}^{-1}.$$
\end{assumption}

Assumption \ref{p1:ass:shrinking_noise} is used only for the explicit
item-count rate. It fails in the persistent subject-effect example
in Proposition \ref{p1:prop:ordinal_noise_floor}.

For the observed combinatorial linearization, write
$$
R_r-\eta_k=H_{rk}+E_{rk},\qquad
H_{rk}=\sum_{b=1}^{G_{k,n}}\sum_{i=1}^{h_{k,b,n}}
C_{k,b,i,\Pi_{r,k,b,n}(i),n}.
$$
Here $G,h,C$ are deterministic given the group, and the arrays
$C_{k,b,ij,n}\in\bbR^T$ refer to the actual ranking model.

\begin{assumption}[Observed permutation modules]\label{p1:ass:signature_linearization}
For fixed $T$ and $n_{\min}\to\infty$, the preceding representation
holds under the actual conditional law, with $h_{k,b,n}\geq2$ and
$\sum_bh_{k,b,n}\to\infty$. The $\Pi_{r,k,b,n}$ are independent uniform
module permutations given $Y_r=k$, and
$$\sum_iC_{k,b,ij,n}=0,\qquad\sum_jC_{k,b,ij,n}=0.$$
\end{assumption}

These permutations are part of a generative representation of the
observed law, not independently introduced calibration permutations.
A common module may contribute to several views. Define
$$
\begin{aligned}
V_{k,n}&=\sum_b\frac{\sum_{i,j}C_{k,b,ij,n}C_{k,b,ij,n}^\top}{h_{k,b,n}-1},\\
c_{k,b,ij,n}(v)&=v^\top V_{k,n}^{-1/2}C_{k,b,ij,n}.
\end{aligned}
$$
and put
$$
\begin{aligned}
\Lambda_{k,n}(v,\epsilon)&=
\sum_b\frac1{h_{k,b,n}}\sum_{i,j}
c_{k,b,ij,n}(v)^2\mathbf1\{|c_{k,b,ij,n}(v)|>\epsilon\},\\
\kappa_{k,n}&=\bbE(\|V_{k,n}^{-1/2}E_{rk}\|_2^2\mid Y_r=k).
\end{aligned}
$$
The inverse is used only under the following nondegeneracy condition.

\begin{assumption}[Gaussian regularity]\label{p1:ass:signature_gaussian}
Under Assumption \ref{p1:ass:signature_linearization}, $V_{k,n}$ is positive
definite. For every active $k$, fixed unit $v\in\bbR^T$, and $\epsilon>0$,
$$n_t(V_{k,n})_{tt}\to\tau_{kt}^2\in(0,\infty),\qquad
\Lambda_{k,n}(v,\epsilon)\to0,\qquad\kappa_{k,n}\to0.$$
\end{assumption}

These are sufficient conditions for an actual-law CLT. They do not
follow from the auxiliary permutation calibration. They are
patternwise, and no uniform Gaussian approximation over growing
dimension or a growing set of groups is asserted.

\subsection{Ordinal information and reference geometry}
\label{p1:subsec:ordinal_information}

\begin{lemma}[Rank invariance and the Spearman identity]
\label{p1:lem:ordinal_signature_identity}
For complete ranks, \eqref{p1:eq:tensor_signature} is well defined and
$$
R^t=\frac{1+\rho_{\mathrm S}(a^t,X^t)}2
=1-\frac{3\|a^t-X^t\|_2^2}{n_t(n_t^2-1)}.
$$
It is unchanged by a strictly increasing transformation applied to all
entries of a raw profile before ranking. Its denominator depends only
on $n_t$.
\end{lemma}

The invariance concerns a common transformation within a profile.
Different item-specific transformations need not preserve its ordering.
The lemma identifies the association statistic exactly; novelty must
therefore lie in the multiview estimator and its recovery theory.

\begin{proposition}[Limits of reference means]
\label{p1:prop:reference_mean_collision}
Different ordinal laws need not have different mean signatures.
For $n=3$, let $P$ assign equal probabilities to $(1,2,3)$ and
$(3,2,1)$, and let $Q$ assign equal probabilities to $(1,3,2)$ and
$(3,1,2)$. Their mean signatures coincide for every rank reference.
Consequently, no reference bank can guarantee mean separation over
the unrestricted class of ranking distributions.
\end{proposition}

The proposition concerns mean-based recovery, not equality of the
entire signature laws. Likewise, iid continuous item values with
different common marginals have the same uniform rank law. Raw
location or scale groups removed by ranking are outside the target.

\begin{theorem}[Multi-reference mean preservation]
\label{p1:thm:multi_reference_identification}
For a finite active set, if every pair has a different mean in at least
one reference coordinate, then the augmented means are injective and
$\Delta_L=\min_{k<\ell}\|\eta_k^{(L)}-\eta_\ell^{(L)}\|_2>0$.
On the deterministic event
$$
\max_{r\leq m}\|R_r^{(L)}-\eta_{Y_r}^{(L)}\|_2<\Delta_L/2,
$$
nearest-true-center assignment recovers every sampled label.
When existing coordinates are retained, an added reference resolves
an old mean collision exactly when its new mean block differs.
\end{theorem}

This is conditional information preservation. It does not establish
the existence of a separating reference bank for arbitrary latent laws.

\begin{proposition}[Overlap with a Mallows Bayes rule]
\label{p1:prop:mallows_overlap}
Consider two conditionally independent multiview Spearman-Mallows
components sharing the reference $a^t$ within each view. Their posterior
log odds have the form
$$
\log\frac{\bbP(Y=1\mid X)}{\bbP(Y=2\mid X)}
=A+\sum_t2D_t(\theta_{1t}-\theta_{2t})R^t,
$$
where $A$ is independent of $X$. If the nonzero coefficients have a
common sign, the nonconstant Bayes rule is a simplex-weighted cutoff,
after orienting the labels. If the threshold falls outside the score
range, the Bayes rule is constant.
\end{proposition}

The result prevents a claim that WRC and ranking-mixture classifiers
are disjoint. It also does not identify a general Mallows mixture with
a fixed-reference one-dimensional classifier.

\subsection{Actual-law concentration and Gaussian limits}
\label{p1:subsec:actual_limits}

\begin{lemma}[Bernstein control from bounded rank modules]
\label{p1:lem:bounded_rank_modules}
Suppose, given $Y=k$, that
$R-\eta_k=\sum_b d_{kb}(\xi_b-p_{kb})$, with independent
$\xi_b\sim\operatorname{Bernoulli}(p_{kb})$ and deterministic vectors
$d_{kb}$. Then Assumption \ref{p1:ass:signature_concentration} holds with
$$
A_kA_k^\top=\sum_b p_{kb}(1-p_{kb})d_{kb}d_{kb}^\top,\qquad
b_k=\frac13\max_b\|d_{kb}\|_1.
$$
\end{lemma}

Shared modules create cross-view covariance. The representation in this
lemma is one sufficient primitive model, not a consequence of ordinal
conversion alone.

\begin{lemma}[Uniform-permutation reference]
\label{p1:lem:permutation_reference}
For complete rank profiles of length $n$, an auxiliary uniform pairing
$\Gamma$ satisfies
$$
\bbE_\Gamma R_\Gamma=\frac12,\qquad
\Var_\Gamma(R_\Gamma)=\frac1{4(n-1)},\qquad
2\sqrt{n-1}(R_\Gamma-\tfrac12)
\overset d\longrightarrow\calN(0,1).
$$
These statements concern calibration, not the observed latent law.
\end{lemma}

\begin{theorem}[Actual-law concentration and Gaussian limit]
\label{p1:thm:signature_limit}
Under Assumption \ref{p1:ass:signature_concentration}, for all $u$ and $x>0$,
$$
\bbP\left(
|u^\top(R_r-\eta_k)|>
\sqrt{2x}\|A_k^\top u\|_2+xb_k\|u\|_\infty
\mid Y_r=k\right)\leq2e^{-x}.
$$
For fixed $T$, Assumptions \ref{p1:ass:signature_linearization} and
\ref{p1:ass:signature_gaussian} imply, conditionally on $Y_r=k$,
$$
V_{k,n}^{-1/2}(R_r-\eta_k)\overset d\longrightarrow\calN_T(0,I_T).
$$
If $\Sigma_{k,n}=\Cov(R_r\mid Y_r=k)$, then
$$
\|V_{k,n}^{-1/2}\Sigma_{k,n}V_{k,n}^{-1/2}-I_T\|_{\mathrm{op}}
\leq2\sqrt{\kappa_{k,n}}+\kappa_{k,n}.
$$
\end{theorem}

The first conclusion is nonasymptotic. The Gaussian conclusions use the
observed-law representation and do not supply a density or likelihood
approximation for the discrete signatures. Recovery below uses only
concentration and is valid at any finite $T$.

For a concrete ordinal family, take $n=2B$, $a=(1,\ldots,2B)$, and
pair item $b$ with item $B+b$. Assign ranks $(b,B+b)$ to the pair when
$\xi_b=0$ and $(B+b,b)$ when $\xi_b=1$, independently with probability
$p$ for a swap. Every generated observation is a complete ranking.

\begin{proposition}[A fully verified pair-swap family]
\label{p1:prop:pair_swap_family}
For the preceding family, put $\lambda_B=3B/(4B^2-1)$ and
$C_B=\sum_b\xi_b\sim\operatorname{Binomial}(B,p)$. Then
$$
R=1-\lambda_B C_B,\quad
\eta_p=1-\lambda_B Bp,\quad
\Var_p(R)=\lambda_B^2Bp(1-p).
$$
For $u>0$,
$\bbP_p(|R-\eta_p|\geq u)\leq2\exp\{-2u^2/(B\lambda_B^2)\}$.
Uniformly away from $p=0,1$, the standardized signature is
asymptotically normal and its variance is of order $n^{-1}$.
\end{proposition}

The noise radius $\lambda_B\sqrt{Bx/2}$ is known from the item count;
it does not require latent labels or a fitted Gaussian density.

\begin{proposition}[Observed block-permutation alternatives]
\label{p1:prop:block_permutation_v2}
For $B\geq2$, take the identity reference of length $2B$ and independent
uniform $\sigma,\tau\in\Sym_B$. Define the aligned profile by
$X_i=\sigma(i)$, $X_{B+i}=B+\tau(i)$, and the reversed profile by
$X_i=B+\sigma(i)$, $X_{B+i}=\tau(i)$.
Their signature means and common variance are
$$
\eta_{+,B}=\frac{7B^2-1}{2(4B^2-1)},\qquad
\eta_{-,B}=\frac{B^2-1}{2(4B^2-1)},\qquad
\sigma_B^2=\frac{(B-1)(B+1)^2}{8(4B^2-1)^2}.
$$
Under either actual law,
$(R-\eta_{\pm,B})/\sigma_B\overset d\longrightarrow\calN(0,1)$.
In particular, $\eta_{+,B}-\eta_{-,B}\to3/4$ and
$2B\sigma_B^2\to1/64$.
\end{proposition}

This is an exact observed-law permutation representation with zero
remainder. It verifies the mean-shift plus asymptotically Gaussian
noise structure for nonuniform full-rank alternatives, rather than
transferring a null permutation CLT to an unspecified latent law.
Fixed collections of independent such views give the corresponding
multivariate limit. The model does not represent arbitrary rank laws.

\begin{lemma}[Multiview pair-swap noise envelopes]
\label{p1:lem:pair_swap_directional_envelopes}
Suppose every view has $n_t=2B$ and, conditionally on each group,
all pair swaps are independent across pairs and views, with arbitrary
group-specific probabilities $p_{kt}$. For every simplex direction,
valid group-uniform projected Bernstein envelopes are
$$
\bar v(w)=\frac{\lambda_B^2B}{4}\|w\|_2^2,\qquad
\bar c(w)=\frac{\lambda_B}{3}\|w\|_\infty.
$$
Alternatively, the same $\bar v(w)$ and $\bar c(w)=0$ form
a valid sub-Gaussian envelope. Thus a sharper valid radius is
$$r^{\rm H}(w,x)=\lambda_B\sqrt{Bx/2}\,\|w\|_2.$$
All envelopes depend only on $B,w$. At URC,
$r^{\rm H}(w,x)$ is of order $\sqrt{x/(BT)}$.
\end{lemma}

The sub-Gaussian alternative follows from the bounded-summand
inequality of \citet{Hoeffding1963Bounded}. It is sharper than the
displayed Bernstein envelope in this particular model; it is not
available merely because general signatures lie in $[0,1]^T$.
These envelopes make direction-specific calibration implementable
without latent labels in this independent-view family. Their use
for shared swaps or a common subject effect would require a
different, verified envelope.

\subsection{Separating projections and variance geometry}
\label{p1:subsec:projection_geometry}

\begin{theorem}[Generic and quantitative simplex separation]
\label{p1:thm:separating_weight_exists}
Suppose $K_0\geq2$ and $\Delta_\infty>0$. Almost every simplex weight
separates the means, and $\gamma_\star=\max_w\gamma(w)>0$.
If $T\geq2$ and $W$ is uniform on $\Delta_T^{\mathrm w}$, then for $g>0$,
\begin{equation}
\bbP\{\gamma(W)\leq g\}
\leq q(g):=\min\left\{1,\frac{4(T-1)\binom{K_0}{2}g}{\Delta_\infty}\right\}.
\label{p1:eq:simplex_small_ball_new}
\end{equation}
For $N$ independent such directions, the probability that all have gap
at most $g$ is at most $q(g)^N$. Also,
$$\gamma(w)\leq\frac1{K_0-1}.$$
\end{theorem}

The quantitative constant depends on dimension and mean geometry.
There is no dimension-free restart or sample-size conclusion.
The final packing bound limits the largest attainable raw margin
when there are many groups. It supplies no lower bound on that margin
and does not identify $K_0$ from a gap loss. In particular, increasing
the number of search directions cannot compensate for a population
margin that is too small relative to the noise.

\begin{proposition}[Recovery using a bank of projections]
\label{p1:prop:projection_bank}
For a fixed bank $w^1,\ldots,w^L$ in the simplex, put
$$
Z_r=( (w^1)^\top R_r,\ldots,(w^L)^\top R_r)^\top,\qquad
\mu_k=( (w^1)^\top\eta_k,\ldots,(w^L)^\top\eta_k)^\top,
$$
and $\Delta_{\rm bank}=\min_{k<\ell}\|\mu_k-\mu_\ell\|_\infty$.
On $\mathcal E_\rho$ of \eqref{p1:eq:uniform_radius_new}, if
$\Delta_{\rm bank}>2\rho$, nearest-true-center
assignment in the sup norm recovers every sampled label. If every
pair admits a direction with projected mean gap at least $h>0$,
a bank of at most $\binom{K_0}{2}$ such directions has
$\Delta_{\rm bank}\geq h$.
\end{proposition}

This is an oracle assignment result in several projection coordinates;
it does not prove unknown-center or unknown-count recovery for a new
multivariate estimator. It explains why retaining several directions
can have a larger robust margin than selecting only one. For example,
the four means $(.2,.2),(.2,.8),(.8,.2),(.8,.8)$ collide along each
coordinate separately and along URC. The single weight $(1/3,2/3)$
separates them with minimum gap $.2$, whereas the two-coordinate
bank has sup-norm separation $.6$. Distinct coordinate values in
each direction are unnecessary for bank-based assignment.

\begin{theorem}[Oracle weighted-cutoff recovery]
\label{p1:thm:oracle_wrc_recovery}
For $K_0\geq2$, under Assumption \ref{p1:ass:signature_concentration}, let $w$ be fixed
with $\gamma(w)>0$ and use the adjacent population midpoint cutoffs.
With $\mathcal I(w)$ defined in \eqref{p1:eq:projected_information},
$$
\bbP\{G_{w,c}(R_r)\neq Y_r\text{ after label alignment}\}
\leq2e^{-\mathcal I(w)},\qquad
\bbP\{\text{any of $m$ errors}\}\leq2me^{-\mathcal I(w)}.
$$
Thus $\mathcal I(w)-\log m\to\infty$ is sufficient for exact recovery.
\end{theorem}

URC obeys the same statement when its projected means differ. A
uniform-weight mean collision is not repaired by moving cutoffs alone.
The bound does not require a separate anti-concentration assumption at
the population cutoffs.
If $\kappa(w)=\gamma(w)/\sqrt{v(w)}$ is already a standardized gap,
the quadratic exponent is of order $\kappa(w)^2$.
The item count enters through $v(w)$; multiplying this exponent
by another item count would generally double-count that scaling.

\begin{corollary}[Accumulation across informative views]
\label{p1:cor:balanced_wrc_rate}
If $v(w)\leq C/(n_{\min}T)$,
$\max_kb_k\leq C/n_{\min}$,
$\|w\|_\infty\leq C/T$, and $\gamma(w)>0$, then
$$n_{\min}T\gamma(w)^2/\log m\longrightarrow\infty$$
is sufficient for oracle exact recovery.
\end{corollary}

The projected covariance condition in this corollary is substantive.
It does not follow from increasing $T$ or from coordinate variance
bounds alone.

\begin{proposition}[Variance-optimal direction in a Gaussian comparison]
\label{p1:prop:variance_geometry}
In the two-group Gaussian comparison experiment with equal priors and
common positive-definite covariance $\Sigma$, put $\delta=\eta_1-\eta_2$.
The midpoint error for direction $w$, with standard normal CDF $\Phi$, is
$$
\Phi\left(-\frac{|w^\top\delta|}{2\sqrt{w^\top\Sigma w}}\right).
$$
The unrestricted optimal squared ratio is $\delta^\top\Sigma^{-1}\delta$.
It is attained by a simplex direction if, after orientation,
$\Sigma^{-1}\delta\geq0$. For
$\Sigma=\sigma_c^2\mathbf1\mathbf1^\top+\sigma_i^2n^{-1}I_T$ and
$\delta=d\mathbf1$, every simplex weight has variance at least
$\sigma_c^2+\sigma_i^2/(nT)$, with equality at URC.
\end{proposition}

This is an exact calculation in the Gaussian comparison, not a
parametric efficiency claim for all ordinal laws. Proposition
\ref{p1:prop:ordinal_noise_floor} gives an actual ordinal counterpart
to persistent common noise.

\subsection{Group count and joint partition recovery}
\label{p1:subsec:joint_gap_recovery}

\begin{proposition}[What a soft gap criterion does not identify]
\label{p1:prop:soft_gap_limit}
For ordered centers in $[0,1]$ with adjacent raw gaps $g_j$ and
$\beta>0$,
$$
\log\sum_{j=1}^{K-1}e^{-\beta g_j}
\geq\log(K-1)-\frac{\beta}{K-1}\qquad(K\geq2).
$$
The increasing lower envelope does not establish monotonicity of
optimized standardized-gap losses across different partitions, or a
unique generating group count. Without component restrictions, even
the full observed distribution need not identify a latent count.
\end{proposition}

In particular, the term $\log\sum_j e^{-8d_j}$ cannot alone justify a
population order-consistency theorem. The certified count instead
uses a noise scale specified independently of a proposed partition.

For $x>0$, define
\begin{equation}
\rho_m(x)=\max_{k,t}\{\sqrt{2x}\,a_{kt}+xb_k\},\qquad
\mathcal E_\rho=\{\max_{r,t}|R_r^t-\eta_{Y_r}^t|\leq\rho\}.
\label{p1:eq:uniform_radius_new}
\end{equation}

\begin{theorem}[Noise-calibrated count and exact partition]
\label{p1:thm:wrc_exact_recovery}
On $\mathcal E_\rho$, every candidate in the simplex satisfies
$\widehat K_\rho(w)\leq K_{\mathrm{obs}}$, the number of sampled groups.
If every active group is sampled and some candidate has
$\gamma(w)>5\rho$, then \eqref{p1:eq:certified_gap_count}--%
\eqref{p1:eq:wrc_cutoff_rule} give
$$
\widehat K=K_0,\qquad
\exists\tau\in\Sym_{K_0}:\ \tau(\widehat Y_r)=Y_r
\quad\text{for all }r\leq m.
$$
Under Assumptions \ref{p1:ass:sampling_abundance} and
\ref{p1:ass:signature_concentration}, take $\rho=\rho_m(x)$.
If $T\geq2$, the pool contains $N$ independent uniform-simplex draws,
and $g>5\rho$, the failure probability is at most
\begin{equation}
2mTe^{-x}+K_0e^{-m\pi_{\min}}+q(g)^N.
\label{p1:eq:joint_gap_error_new}
\end{equation}
Additional data-dependent candidates do not invalidate this bound.
For $K_0=1$ the projection-coverage term is unnecessary; for $T=1$ it
is unnecessary whenever $\gamma(1)>5\rho$.
\end{theorem}

The theorem identifies a partition, not a unique weight vector.
It proves a property of the displayed estimator, without postulating
that a soft-loss minimizer is already oracle-equivalent.

\begin{corollary}[An explicit, conservative random-pool regime]
\label{p1:cor:random_pool_rate}
Suppose Assumption \ref{p1:ass:shrinking_noise} holds, $T,K_0\geq2$, and
$x_m\to\infty$ satisfies $mTe^{-x_m}\to0$. If
$$
\frac{TK_0^2}{\Delta_\infty}
\left(\sqrt{\frac{x_m}{n_{\min}}}+\frac{x_m}{n_{\min}}\right)\to0,
\quad N\to\infty,\quad
m\pi_{\min}-\log K_0\to\infty,
$$
then the joint recovery probability tends to one.
\end{corollary}

This is a sufficient bound for a universal random pool, not a minimax
requirement. Coordinate or pilot directions can be substantially better
when the signal is concentrated in a small set of views.

\begin{theorem}[Direction-specific noise calibration]
\label{p1:thm:directional_calibration}
Under Assumptions \ref{p1:ass:sampling_abundance} and
\ref{p1:ass:signature_concentration}, let $\mathcal W$ be a fixed
nonempty bank of size $J$, with valid envelopes
$\bar v_j\geq v(w^j)$, $\bar c_j\geq\beta(w^j)$ and radii $r_j(x)$
defined above. Apply the extended gap rule with threshold
$\vartheta r_j(x)$, where $\vartheta\geq2$. With probability at least
$1-2mJe^{-x}$, no candidate splits a sampled group. If, in addition,
$$
\mathcal G_x=
\{\exists j:\gamma(w^j)>(\vartheta+2)r_j(x)\}
$$
holds and every active group is sampled, the selected count and
partition are correct. For a random bank independent of the recovery
data, with conditionally valid envelopes, the failure probability is
at most
$$
2mJe^{-x}+K_0e^{-m\pi_{\min}}+\bbP(\mathcal G_x^c).
$$
The same statement holds conditionally on independent pilot data.
If the conditional envelopes fail with probability at most
$\delta_{\rm cal}$, add $\delta_{\rm cal}$ to the bound.
The result also applies to any other radii satisfying the same
conditional $2e^{-x}$ scalar tail bound.
\end{theorem}

The proof uses one simultaneous projection event. On that event,
each true group has diameter at most $2r_j$ in coordinate $j$, so
the calibrated cuts cannot split it. One adequate coordinate separates
all groups, and count maximization then recovers the partition.
The choice $\vartheta=2$ is valid because cuts use a strict inequality;
it improves the sufficient gap from $5r_j$ to $4r_j$, without
making a worst-case radius exact at a finite sample size.

For $N$ iid proposals with the same deterministic radius map $r_x(w)$,
define their probability of an adequate direction by
$$
p_x=\bbP_W\{\gamma(W)>(\vartheta+2)r_x(W)\}.
$$
If these proposals are included in the bank, their coverage failure is
$(1-p_x)^N\leq e^{-Np_x}$. The probability must be bounded
conditionally when the proposal distribution depends on an independent
pilot. This identifies a useful target for candidate learning: the
probability of a noise-resolved margin, rather than mere nonzero
population separation. For a general random bank, the theorem keeps
$\bbP(\mathcal G_x^c)$ explicit and assumes no independence among its
directions.

\begin{corollary}[Growing group count with an attainable margin]
\label{p1:cor:growing_group_count}
Suppose Assumptions \ref{p1:ass:sampling_abundance},
\ref{p1:ass:signature_concentration}, and
\ref{p1:ass:shrinking_noise} hold, and $K_0=K_m\geq2$.
Let $\alpha_m\to0$ and $x_m=\log(2mT/\alpha_m)$.
Suppose, for a constant $c_0>0$, the candidate bank satisfies
$$
\bbP\left\{\max_{w\in\mathcal W}\gamma(w)
             \geq\frac{c_0}{K_m-1}\right\}\longrightarrow1.
$$
If
$$
\frac{n_{\min}}{(K_m-1)^2x_m}\longrightarrow\infty,
\qquad m\pi_{\min}-\log K_m\longrightarrow\infty,
$$
then the original uniform-radius estimator jointly recovers the count
and partition with probability tending to one.
\end{corollary}

The prescribed lower margin is an additional geometric and search
condition; it is not implied by the packing upper bound. For balanced
groups with $\pi_{\min}\geq c/K_m$, a convenient sufficient occupancy
condition is $m/(K_m\log K_m)\to\infty$. The corollary describes a
bank attaining a margin of order $K_m^{-1}$. By comparison, the
universal random-pool condition in Corollary
\ref{p1:cor:random_pool_rate} is ensured by
$$
\frac{n_{\min}\Delta_\infty^2}
     {T^2K_m^4x_m}\longrightarrow\infty
$$
together with its occupancy and restart conditions. This stronger
sufficient rate reflects the conservative simplex small-ball bound;
it is not a universal lower bound on the information required.

\begin{corollary}[Joint recovery across informative views]
\label{p1:cor:directional_view_recovery}
Consider the independent-bank estimator of Theorem
\ref{p1:thm:directional_calibration}, with fixed $\vartheta\geq2$,
$K_0=K_m\geq2$, $J=J_m$, $\alpha_m\to0$, and
$x_m=\log(2mJ_m/\alpha_m)$. Suppose the bank contains, with probability
tending to one, a direction satisfying
$$
\gamma(w)\geq\frac{c_0}{K_m-1},\qquad
\bar v(w)\leq\frac{C}{n_{\min}T},\qquad
\bar c(w)\leq\frac{C}{n_{\min}T}.
$$
Then joint count and partition recovery follows if
$$
\frac{n_{\min}T}{(K_m-1)^2x_m}\to\infty,
\qquad m\pi_{\min}-\log K_m\to\infty.
$$
\end{corollary}

The view-accumulation condition concerns a validated projected
variance and tail scale. Independent informative views can satisfy
it, whereas replicated views need not. A persistent subject effect
can invalidate even coordinate noise shrinkage. For fixed $K_0,T$
and a stable positive gap, increasing the item counts sufficiently
fast relative to $\log m$ yields joint recovery; increasing $m$
alone does not remove overlapping subject laws. For a fixed oracle
classifier with iid subject error probability $p_m$, the probability
of all named labels being correct is $(1-p_m)^m$. Thus vanishing
individual error is weaker than complete recovery, whose relevant
small-error condition is $mp_m\to0$.

\begin{corollary}[Prediction from certified sample cutoffs]
\label{p1:cor:certified_prediction}
On the successful training event of Theorem
\ref{p1:thm:wrc_exact_recovery}, suppose $K_0\geq2$ and
$\widehat g>4\rho$. For an independent test subject, put
$a=(\widehat g-4\rho)/2$. Conditionally on the training data,
$$
\bbP(\text{test error}\mid\mathcal D)
\leq2\exp\left\{-\frac{a^2}{2\{v(\widehat w)+\beta(\widehat w)a\}}\right\},
$$
where $\beta(w)=\max_kb_k\|w\|_\infty$.
If a candidate has population gap greater than $6\rho$, the specified
gap tie-break guarantees $\widehat g>4\rho$ on that event.
\end{corollary}

The statement concerns the successful training event. Unconditional
risk also includes the training failure bound
\eqref{p1:eq:joint_gap_error_new}. No fitted empirical margin is
substituted for an unverified population concentration scale.

For a fitted nonempty group $\widehat C_j$, define its estimated signature
center by $\widehat\eta_j=|\widehat C_j|^{-1}\sum_{r\in\widehat C_j}R_r$.
Let $N_k=\sum_r\mathbf1\{Y_r=k\}$.

\begin{corollary}[Recovery of the signature means]
\label{p1:cor:signature_mean_recovery}
Under Assumptions \ref{p1:ass:sampling_abundance} and
\ref{p1:ass:signature_concentration}, for $x>0$, the probability that the
count and partition are correctly recovered but at least one inequality
$$
|\widehat\eta_{\tau^{-1}(k),t}-\eta_k^t|
\leq\sqrt{\frac{2x}{N_k}}\,a_{kt}+\frac{xb_k}{N_k}
\quad(k=1,\ldots,K_0,\ t=1,\ldots,T)
$$
fails is at most $2K_0Te^{-x}$.
Adding the failure bound \eqref{p1:eq:joint_gap_error_new} gives joint
partition and signature-center recovery. Joint consistency follows
if that failure bound vanishes and, for some $x_m$ with
$K_0Te^{-x_m}\to0$, the maximum displayed radius tends to zero in probability.
\end{corollary}

The intersection with correct recovery is essential: the proof controls
true-label sample means and then transfers them on the recovery event,
without treating selected clusters as an independent sample.

\subsection{Recovery from combinations of coordinates}
\label{p1:subsec:tuple_theory_v2}

For bank coordinate $\ell$, write its distinct population levels as
$a_{\ell1}<\cdots<a_{\ell q_\ell}$, where
$\{a_{\ell j}:j\leq q_\ell\}=\{(w^\ell)^\top\eta_k:k\leq K_0\}$.
Mark a population cut $j$ when
$a_{\ell,j+1}-a_{\ell j}>(\vartheta+2)r_\ell(x)$, and define
$$
\mathcal J_\ell(x)=\{j<q_\ell:
 a_{\ell,j+1}-a_{\ell j}>(\vartheta+2)r_\ell(x)\},\qquad
c_\ell(k)=\sum_{j\in\mathcal J_\ell(x)}
                  \mathbf1\{(w^\ell)^\top\eta_k>a_{\ell j}\}.
$$
The population code $c(k)=(c_1(k),\ldots,c_L(k))$ is used only in
the theorem. The estimator in \eqref{p1:eq:tuple_estimator_v2}
does not calculate it.

\begin{theorem}[Unknown-count recovery by coordinate refinement]
\label{p1:thm:tuple_recovery_v2}
Let the fixed or independently generated bank have valid conditional
scalar tail radii $r_\ell(x)$, as in Theorem
\ref{p1:thm:directional_calibration}. With probability at least
$1-2mLe^{-x}-\delta_{\rm cal}$, every tuple block is a union of sampled true groups.
If $k\mapsto c(k)$ is injective, that event also gives the exact
partition of sampled subjects. Under Assumption
\ref{p1:ass:sampling_abundance}, joint recovery of the active count
and partition fails with probability at most
$$
2mLe^{-x}+K_0e^{-m\pi_{\min}}
       +\bbP\{k\mapsto c(k)\text{ is not injective}\}
       +\delta_{\rm cal}.
$$
Here $\delta_{\rm cal}$ is zero for surely valid envelopes. No
individual coordinate is required to have $K_0$ distinct levels.
\end{theorem}

The resolved population cuts ensure that every pair is separated in
at least one empirical coordinate partition. The no-splitting event
ensures that their common refinement cannot fragment a true group.
This proves unknown-center recovery for the implemented tuple rule,
under a stronger coordinate-level condition than the nearest-center
bank result. The number of groups can be as large as
$\prod_\ell q_\ell$ when the corresponding level combinations occur.

\begin{proposition}[Unknown-center recovery by a calibrated graph]
\label{p1:prop:bank_graph_v2}
Use positive radii and the graph defined in Section
\ref{p1:subsec:tuple_method}. On the simultaneous projection event,
its components are unions of sampled true groups. If
$$
\min_{k<h}\max_{\ell\leq L}
\frac{|(w^\ell)^\top(\eta_k-\eta_h)|}{r_\ell(x)}>4,
$$
the components are exactly the sampled groups. Joint count recovery
then has failure bound $2mLe^{-x}+K_0e^{-m\pi_{\min}}$, plus the
probability that the displayed bank condition or calibration fails.
\end{proposition}

The graph does not require gaps between successive population levels;
the tuple theorem requires its resolved cuts. Both statements are
finite-sample consequences of calibrated geometry, rather than
generic consistency assertions for single linkage on overlapping laws.

\begin{corollary}[A replicated binary-code geometry]
\label{p1:cor:binary_code_v2}
In the independent pair-swap model, partition $T=qh$ views into
$q$ fixed blocks of size $h$, and average within each block.
The active classes have distinct codes $v_k\in\{0,1\}^q$ and swap
probabilities $.2+.6v_{k\ell}$ in block $\ell$.
For $\alpha_m\to0$, set $x_m=\log(2mq/\alpha_m)$ and use
$$r_\ell(x_m)=\lambda_B\sqrt{Bx_m/(2h)}.$$
Tuple count and partition recovery follows if
$$Bh/x_m\to\infty,\qquad
 m\pi_{\min}-\log K_0\to\infty.$$
The blocks are specified independently of the active subset of codes.
\end{corollary}

Each nonconstant block has two means separated by $.6\lambda_BB$.
Consequently no block alone need identify all classes. With
$q$ of order $\log K_0$ and $h=T/q$, the sufficient information
condition is $BT/\{q\log(2mq/\alpha_m)\}\to\infty$.
This example avoids the one-dimensional inverse-$K_0$ margin loss;
it assumes a known grouping of views into code coordinates and does
not establish a general minimax rate over unknown geometries.

\subsection{Ordinal lower bounds and persistent noise}
\label{p1:subsec:ordinal_lower_bounds}

For $p,q\in(0,1)$, write
$\operatorname{kl}(p,q)=p\log(p/q)+(1-p)\log\{(1-p)/(1-q)\}$.

\begin{theorem}[The logarithmic recovery boundary in an ordinal submodel]
\label{p1:thm:ordinal_lower_bound}
In the pair-swap family, take two equally likely groups with
$p_+,p_-\in[\epsilon,1-\epsilon]$, $\delta=|p_+-p_-|>0$, and
$D=B\,\operatorname{kl}(p_+,p_-)$, where $\operatorname{kl}$ is Bernoulli
relative entropy. For any estimator, even knowing the component laws,
$$
\bbP(\text{all labels correct up to permutation})
\leq2\exp\{-(m/4)e^{-D}\}.
$$
Hence $D\leq(1-\varepsilon)\log m$ implies failure of exact recovery
as $m\to\infty$. Conversely, the calibrated gap estimator with $T=1$
recovers the unknown count and partition when
$B\delta^2/\log m\to\infty$. The necessary and sufficient scales agree
up to constants; no sharp universal constant is claimed.
For $T$ conditionally independent copies with the same probabilities,
replace $B$ by $BT$ in the information and use a direction-calibrated
URC rule for the sufficient condition.
\end{theorem}

\begin{proposition}[An exact known-law binary comparison]
\label{p1:prop:exact_bayes_v2}
In the symmetric independent pair-swap model with equal priors, put
$p_-=(1-\delta)/2$, $p_+=(1+\delta)/2$, $0<\delta<1$, and $L=BT$.
The total swap count is sufficient for the full ranking observations.
With independent fair randomization at a tied count, the Bayes
individual error and the partition success of pointwise MAP are
$$
p_e=\bbP_{p_-}\{\operatorname{Bin}(L,p_-)>L/2\}
 +\tfrac12\bbP_{p_-}\{\operatorname{Bin}(L,p_-)=L/2\},\qquad
(1-p_e)^m+p_e^m.
$$
The second expression evaluates this classifier, rather than proving
that it is the optimal estimator of an unlabeled partition.
\end{proposition}

This supplies an exact benchmark for the aligned binary experiment.
A direction selected by maximum raw gap and given population cutoffs
is a different, restricted known-center comparison, not a Bayes upper
envelope over all directions or decision rules.

\begin{proposition}[An ordinal variance floor]
\label{p1:prop:ordinal_noise_floor}
Modify the pair-swap family by drawing a subject effect $U$ and then
independent swaps with probability $U$ in every view, conditionally on $U$.
With equal group priors, let
$U\mid Y=1\sim\operatorname{Unif}[.1,.6]$ and
$U\mid Y=2\sim\operatorname{Unif}[.4,.9]$.
Every classifier has individual error at least $.2$, at any item
counts or number of such views. Even with known component laws, any
partition estimator has success probability at most $2(.8)^m$. For each group,
$$
\Var(R^t\mid Y=k)\longrightarrow
\frac9{16}\Var(U\mid Y=k)>0.
$$
Thus distinct mean signatures do not imply exact recovery when
common subject heterogeneity persists.
\end{proposition}

The Gaussian and ordinal calculations separate reducible within-view
variation from subject-level overlap. Assumptions that force variance
to vanish must be verified against this distinction.

\begin{proposition}[New ordinal information and a lower floor]
\label{p1:prop:anchor_floor_v2}
In Proposition \ref{p1:prop:ordinal_noise_floor}, add a binary rank
view encoding $A\in\{0,1\}$ with
$\bbP(A\ne Y\mid Y,U)=\epsilon\in[0,1/2]$; here labels are $0,1$.
Given $U$, the anchor is independent of the base rank modules.
The oracle observing $(U,A)$ has error $.4\epsilon$.
For $L=BT$ base modules, the fixed rule
$\widehat Y=\mathbf1\{C/L>.4\}$ if $A=1$, and
$\widehat Y=\mathbf1\{C/L>.6\}$ if $A=0$, obeys, for $a>0$,
$$
.4\epsilon\leq\bbP(\widehat Y\ne Y)
\leq .4\epsilon+8a+2e^{-2La^2}.
$$
It is a fixed simplex WRC rule after encoding the anchor score as $1-A$.
For fixed $\epsilon$, its error tends to $.4\epsilon$ as $L\to\infty$.
\end{proposition}

Its finite-sample lower bound follows from the enlarged-information
oracle; it does not contradict the original $.2$ floor.
For $\epsilon=\epsilon_m$, the displayed upper bound also gives
complete named-label recovery if $m\epsilon_m\to0$ and
$BT/(m^2\log m)\to\infty$. This last sufficient rate is deliberately
conservative, not a sharp boundary. The construction uses known
overlap endpoints to design its weights and cutoff, not learned
unlabeled weights.

%------------------------------------------------------------
%Simulation Study
%------------------------------------------------------------

\section{Simulation Studies}
\label{p1:sec:simulation}
\subsection{Design and reporting conventions}
\label{p1:sim:design}

We examine five focused scenarios directly linked to the main recovery
results. They are ordered by their role in the analysis: calibrated
cutoffs, increasing group counts, coordinate collisions, view-dependent
variance, and an additional ordinal anchor. These are illustrative
simulation comparisons selected from a broader exploratory study, not an
exhaustive benchmark or a claim of uniformly favorable performance.
Simulation results are reported only for numerical cells documented in
the available summary. The complete cell-level summary for Scenario~1 is available, but
replication-level outputs and complete cell-level summaries for the other
scenarios were not supplied; consequently no unreported intermediate
points or confidence intervals are reconstructed.

The studies use the indicated subsets of
\[
T\in\{1,4,16,64,256,1024\},\qquad
B\in\{8,32,128,512,2096\},
\]
with scenario-specific exceptions for the adaptive $B$ used in
Scenario~2. In the complete-rank pair-swap design, $B$ counts independent
rank pairs ($n_t=2B$); exact binomial counts are converted to scores by
$R=1-\lambda_B C$. The scenarios state their additional distributional,
weight-bank, threshold, and sample-size choices below; definitions of
the theoretical radii appear in Section~\ref{p1:sec:theory}. The
horizontal value $1$ in the recovery panels is a descriptive
perfect-recovery reference, \emph{not} an asserted theoretical threshold.
An observed recovery proportion of $1$ over 2,000 independent
replications has a single-cell 95\% Wilson lower confidence bound of
approximately $0.9981$.

\subsection{Scenario 1: calibration with a fixed candidate direction}
\label{p1:sim:scenario1}

This experiment isolates the effect of noise calibration while holding
both the candidate direction and the cutoff multiplier fixed. We use
the aligned binary pair-swap model with $K=2$, $m=1000$, $B=32$,
and $\vartheta=3$. Both procedures project onto the same uniform
reference-coordinate direction,
$w^{\mathrm U}=T^{-1}(1,\ldots,1)^\top$.
The baseline uses the coordinatewise Hoeffding radius, while the
alternative applies a direction-specific Bernstein radius to the
projected score. Thus, neither the weight nor the clustering rule
changes across the two procedures.

Table~\ref{p1:tab:scenario1_calibration} reports the implemented
noise radii and joint exact-recovery proportions from 2,000
Monte Carlo repetitions per cell. Both calibration rules fail to
recover the full partition and group count in all observed repetitions
at $T=4$ and $T=16$. Although the projected Bernstein radius is
smaller in those cells, it remains too conservative to produce
exact recovery. From $T=64$ onward, the projected Bernstein
rule succeeds in all 2,000 repetitions, whereas the
coordinate-radius rule records no exact recoveries. In particular,
at $T=64$ the two radii are approximately $0.4052$ and $0.0464$,
respectively. This contrast illustrates how accounting for
independent-view averaging can substantially improve the calibration
of a fixed projection. The observed change between $T=16$ and
$T=64$ is a finite-grid phenomenon, not an established universal
recovery threshold; likewise, an observed proportion of one does
not imply zero population error probability.

\begin{table}[!htbp]
\centering
\caption{Scenario 1: implemented noise radii and joint exact-recovery
proportions under two calibrations of the same uniform projection.
The aligned binary pair-swap design fixes $K=2$, $m=1000$, $B=32$,
and $\vartheta=3$. Each cell has 2,000 Monte Carlo repetitions.
The coordinate calibration uses a Hoeffding bound; the projected
calibration uses a direction-specific Bernstein bound.}
\label{p1:tab:scenario1_calibration}
\begin{tabular}{rcccc}
\toprule
 & \multicolumn{2}{c}{Noise radius} & \multicolumn{2}{c}{Joint exact recovery} \\
\cmidrule(lr){2-3}\cmidrule(lr){4-5}
$T$ & Coordinate & Projected Bernstein & Coordinate & Projected Bernstein \\
\midrule
4 & 0.3739 & 0.2069 & 0.0000 & 0.0000 \\
16 & 0.3898 & 0.0964 & 0.0000 & 0.0000 \\
64 & 0.4052 & 0.0464 & 0.0000 & 1.0000 \\
256 & 0.4199 & 0.0228 & 0.0000 & 1.0000 \\
1024 & 0.4342 & 0.0113 & 0.0000 & 1.0000 \\
\bottomrule
\end{tabular}
\end{table}

\subsection{Scenario 2: increasing views with a growing-group design}
\label{p1:sim:scenario2}

\begin{figure}[!htbp]
\centering
\includegraphics[width=.86\linewidth]{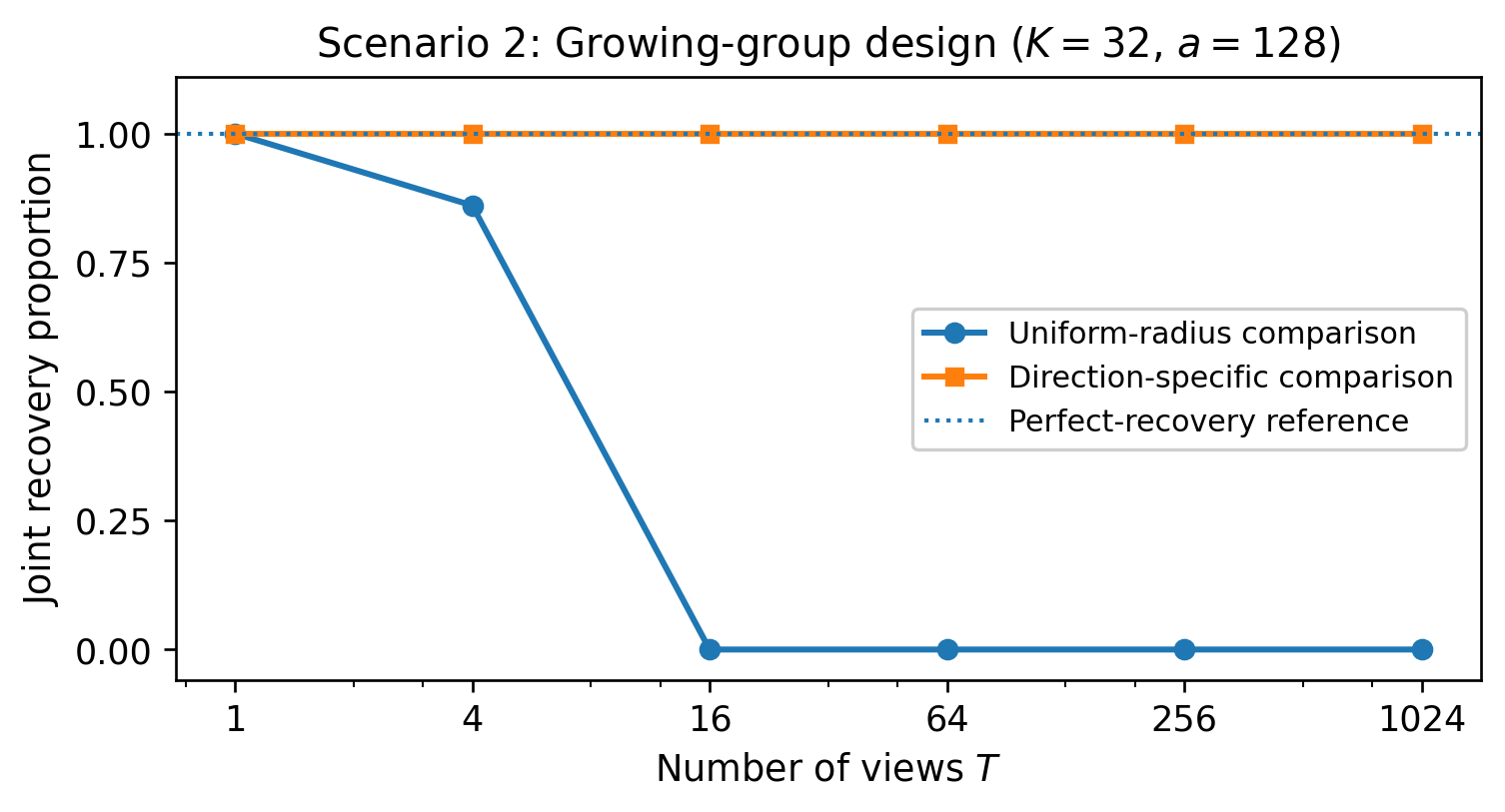}
\caption{Scenario 2: observed joint-recovery proportions for
$K=32$, $a=128$ (2,000 repetitions per reported point). The dotted
line is the perfect-recovery reference. Candidate banks and cutoff
multipliers differ between the two curves, as specified in
Table~\ref{p1:tab:scenarios23}; $B$ changes with $T$.
Lines connect reported cell summaries only.}
\label{p1:fig:scenario2}
\end{figure}

Set $K=32$, and let the group-specific swap probabilities be equally
spaced over $[0.15,0.85]$ with alternating view orientations. A specified
parity-average direction has population gap of order $(K-1)^{-1}$.
The study sets
\[
m=\lceil 30K\log(K+1)\rceil,\qquad
B=\max\!\left\{8,\left\lceil
\frac{a(K-1)^2\log(m)}{T}\right\rceil\right\},
\qquad a=128.
\]
Figure~\ref{p1:fig:scenario2} and Table~\ref{p1:tab:scenarios23}
show the reported path over six values of $T$. The specified
direction-specific method achieves complete observed recovery over this
path while the reported uniform-radius alternative deteriorates as
$T$ increases. \emph{These columns also differ in cutoff multiplier
and candidate bank}; thus this scenario does not separately identify
the benefit of calibration. That more controlled comparison is
Scenario~1. The result is compatible with Corollary
\ref{p1:cor:directional_view_recovery} under its stated conditions,
but does not establish a universal finite-sample cutoff.

\subsection{Scenario 3: recovery through coordinate refinement}
\label{p1:sim:scenario3}

\begin{figure}[!htbp]
\centering
\includegraphics[width=.86\linewidth]{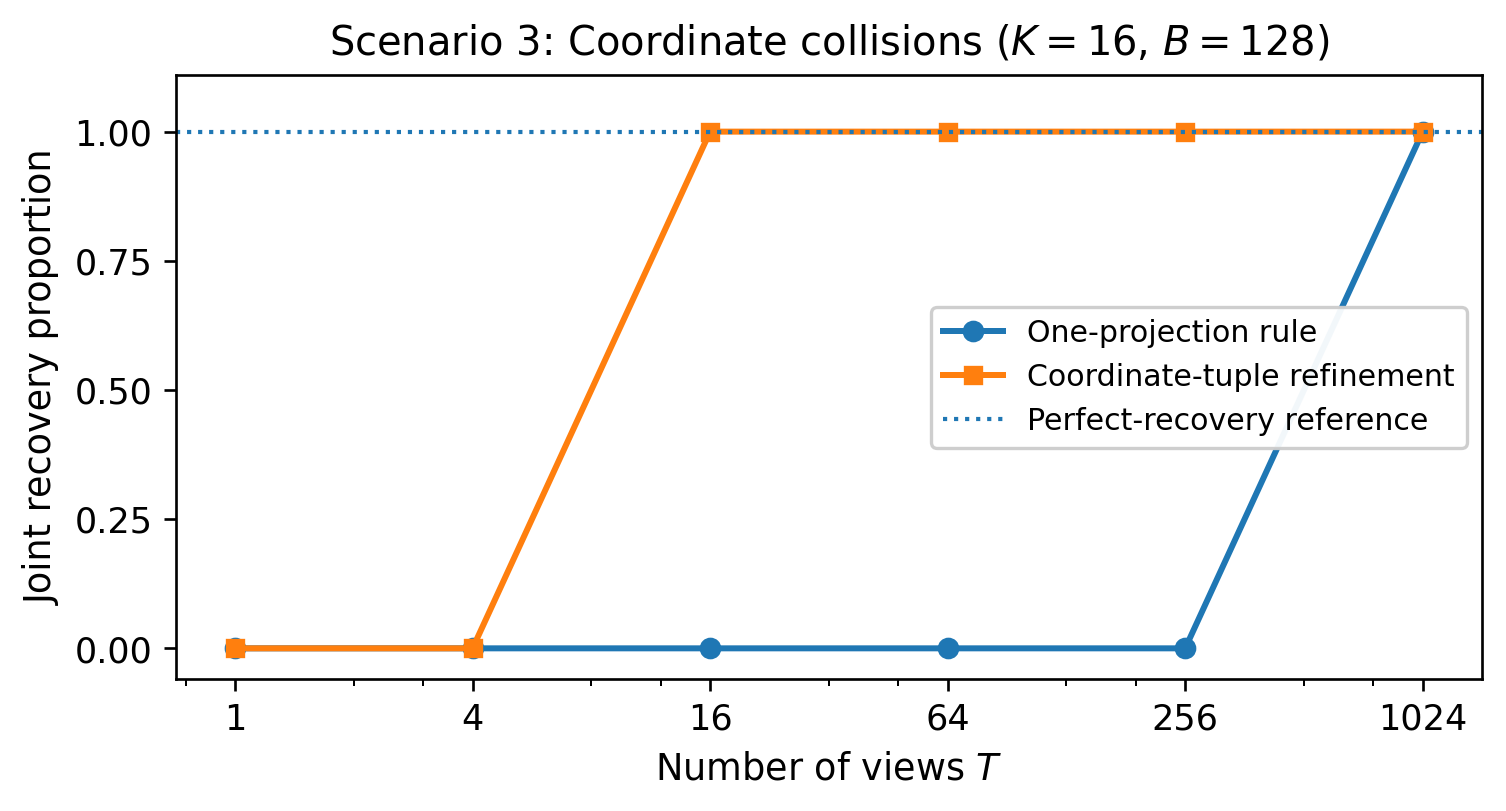}
\caption{Scenario 3: observed joint recovery with $K=16$,
$B=128$ and 2,000 repetitions per point. The dotted horizontal line
denotes the descriptive perfect-recovery reference. The projection
bank and cutoff follow Table~\ref{p1:tab:scenarios23}.}
\label{p1:fig:scenario3}
\end{figure}

Select $K=16$ distinct codes from a four-bit view layout. Each bit-view
has swap probability $0.2$ or $0.8$, with $B=128$ independent rank pairs
per view. Every original coordinate has at most two population levels
and therefore has collisions across the sixteen groups. The feasible
tuple procedure estimates coordinate-level groups before clustering
subjects by tuples of coordinate labels; it is not supplied the true
centers or $K$. In the documented grid, joint recovery becomes complete
from $T=16$ onward, whereas the augmented single-projection rule first
reaches complete recovery at $T=1024$ (Figure~\ref{p1:fig:scenario3}).
This illustrates the resolved-level condition of Theorem
\ref{p1:thm:tuple_recovery_v2}; it does not imply that tuple refinement
universally dominates optimized one-dimensional projections.

\begin{table}[!htbp]
\centering
\caption{Reported joint-recovery proportions for Scenarios 2 and 3
(2,000 repetitions per cell). Scenario 2 uses $K=32$, $a=128$ and
adaptive $B$; Scenario 3 uses $K=16$ and fixed $B=128$.
The direction-specific and projection columns use the augmented
candidate bank, Hoeffding radii and $\vartheta=2$; the Scenario 2
uniform-radius column uses $\vartheta=3$ and its original bank.
The differing settings preclude attributing the Scenario 2 contrast
solely to calibration.}
\label{p1:tab:scenarios23}
\begin{tabular}{rrrrrr}
\toprule
& \multicolumn{3}{c}{Scenario 2} & \multicolumn{2}{c}{Scenario 3}\\
$T$ & $B$ & Uniform & Direction-specific & Projection & Tuple\\
\midrule
1 & 998678 & 1.0000 & 1.0000 & 0.0000 & 0.0000\\
4 & 249670 & .8605 & 1.0000 & 0.0000 & 0.0000\\
16 & 62418 & 0.0000 & 1.0000 & 0.0000 & 1.0000\\
64 & 15605 & 0.0000 & 1.0000 & 0.0000 & 1.0000\\
256 & 3902 & 0.0000 & 1.0000 & 0.0000 & 1.0000\\
1024 & 976 & 0.0000 & 1.0000 & 1.0000 & 1.0000\\
\bottomrule
\end{tabular}
\end{table}

\subsection{Scenario 4: independent versus replicated views}
\label{p1:sim:scenario4}

\begin{figure}[!htbp]
\centering
\includegraphics[width=.86\linewidth]{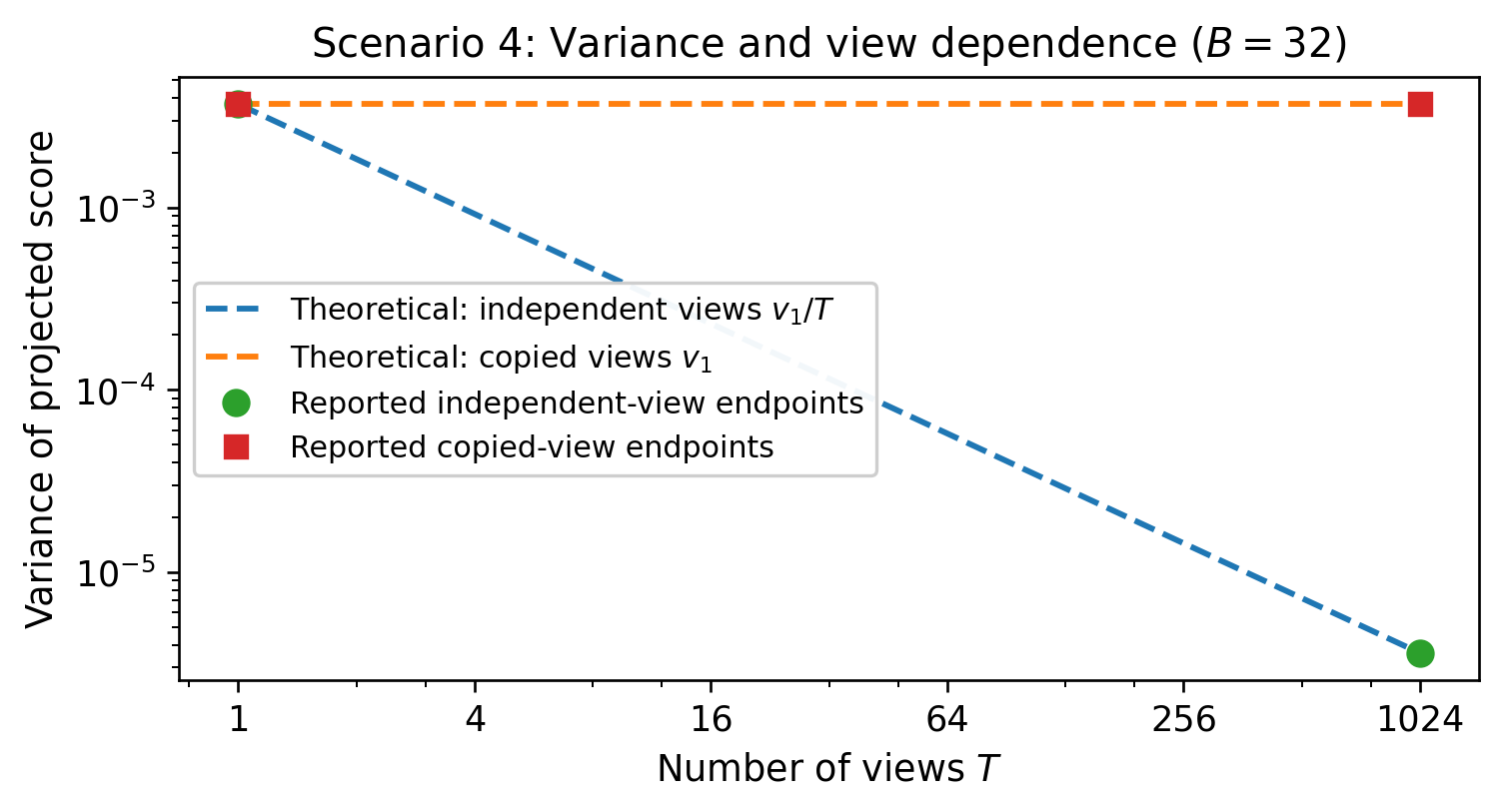}
\caption{Scenario 4: cross-view variance benchmarks. Dashed curves
represent $v_1/T$ for independent views and $v_1$ for copied views.
Only endpoint markers represent reported Monte Carlo variance values;
$v_1$ is normalized to the reported $T=1$ value.}
\label{p1:fig:scenario4}
\end{figure}

With $B=32$, the reported variance of a uniform score projection
under independent views is approximately $0.003693$ at $T=1$ and
$0.000003606$ at $T=1024$. Perfectly copied views retain variance
approximately $0.003693$. Figure~\ref{p1:fig:scenario4} compares
these observed endpoints with the theoretical variance $v_1/T$
for independent identically distributed views and the constant $v_1$
for perfectly copied views, using $v_1=0.003693$. The intermediate
segments are \emph{theoretical reference curves, not simulated points}.
This example makes the dependence assumptions behind projected
calibration visible.

\subsection{Scenario 5: reducing an irreducible error floor with an anchor}
\label{p1:sim:scenario5}

\begin{figure}[!htbp]
\centering
\includegraphics[width=.86\linewidth]{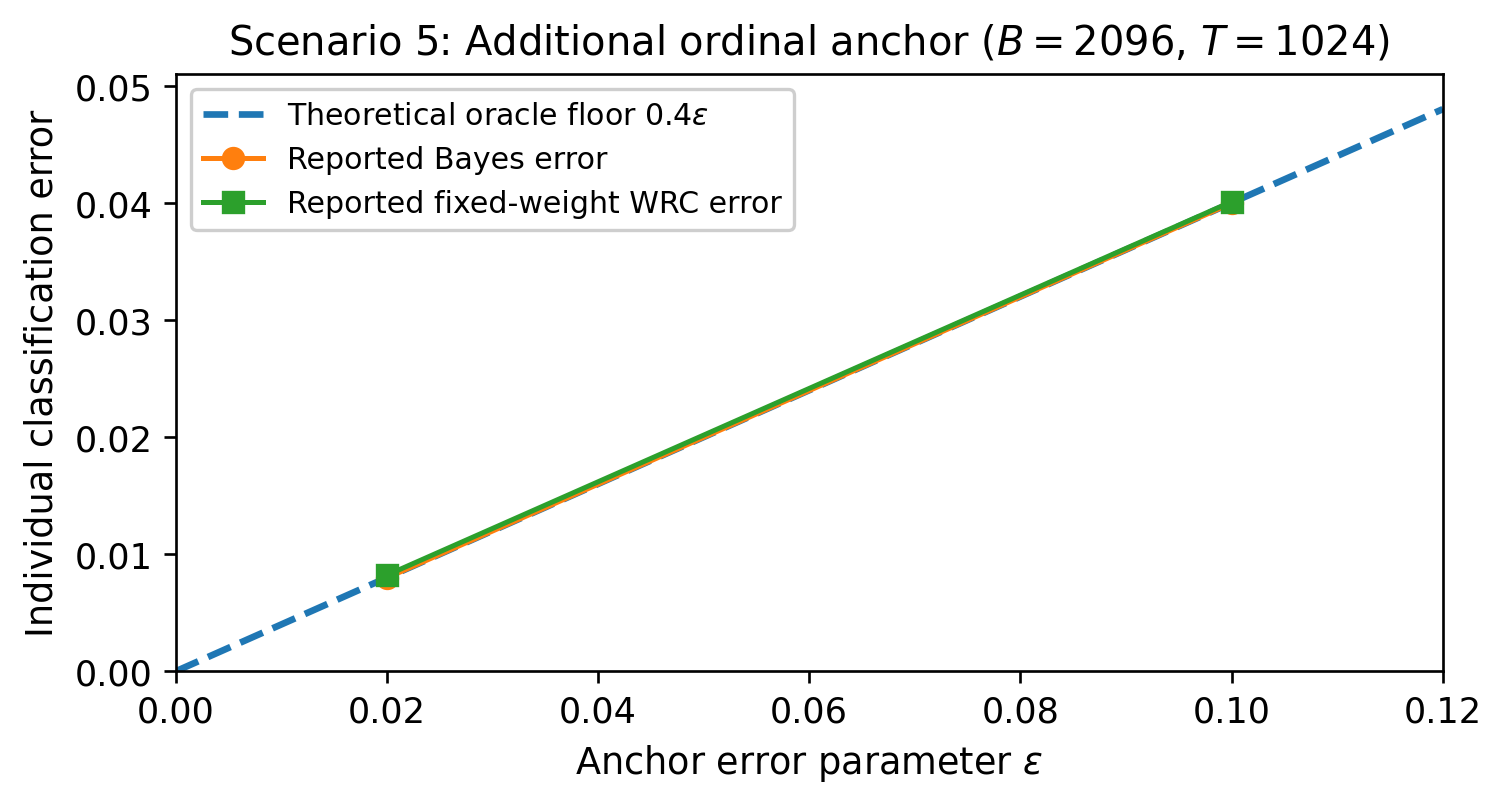}
\caption{Scenario 5: individual classification errors for the
fixed-weight WRC and enlarged-information Bayes procedures against the
proven oracle floor $0.4\epsilon$. The line is theoretical, not fitted
to the two observed cells ($B=2096,T=1024$); segments between
measured points are visual guides only.}
\label{p1:fig:scenario5}
\end{figure}

In the shared-subject-effect construction, the rank-only model has an
individual-error floor of $0.2$. The additional binary ordinal anchor
of Proposition~\ref{p1:prop:anchor_floor_v2} gives an
enlarged-information Bayes floor of $0.4\epsilon$. At $B=2096$ and
$T=1024$, the reported anchor-model Bayes errors are $0.007976$ and
$0.039967$ for $\epsilon=0.02$ and $0.10$, respectively; the
fixed-weight WRC errors are $0.008210$ and $0.040096$.
Figure~\ref{p1:fig:scenario5} compares these two measured cells with
the proven $0.4\epsilon$ line. The weights and cutoff use known
overlap endpoints, so this does not assess unsupervised weight learning.

\subsection{Interpretation and limitations}
The five scenarios test distinct sufficient-condition mechanisms, not
an end-to-end fitted model or general empirical dominance. In
particular, the directional and uniform comparisons in Scenario 2 are
not a single-factor ablation, and the constructed anchor in Scenario 5
uses known distributional parameters. Negative outcomes and unreported
settings of the broader exploratory campaign have not been converted
into successes; the displayed evidence is a selected illustration and
cannot be treated as a comprehensive method comparison. No fitted
Mallows-mixture benchmark, neural training study, or bootstrap-weight
learning evaluation is claimed here. The simulation scenarios are interpreted only within their
specified data-generating models and calibration settings.

\section{Conclusion}
\label{p1:sec:conclusion}

We studied exact recovery of the number of latent groups and the
partition of subjects from reference-based ordinal rearrangement
signatures. The representation is deliberately rank-invariant; on
complete rankings its individual coordinates are affine Spearman
correlations. The statistical contribution is therefore not a new rank
correlation, but a description of what multiview reference signatures
retain and how calibrated noise controls exact recovery.

For the explicitly defined finite candidate bank, the analysis
separates three obstacles: nonidentifiable signature means, random
fluctuation of projected observations, and failure to sample a
sufficiently separating projection. Actual-law concentration bounds
yield nonasymptotic recovery guarantees for both the group count and
subject partition. The finite-dimensional Gaussian limit has a
separate role and is not invoked as a growing-dimensional likelihood
approximation. An ordinal information lower bound explains the
logarithmic dependence on the subject count within a specified
submodel, without asserting global minimax optimality.

Direction-specific calibration can exploit variance reduction across
informative independent views, whereas correlated or copied views do
not necessarily provide additional information. When many groups must
be distinguished, a single bounded projection faces a diminishing
minimum-gap constraint; coordinatewise partitions and their common
refinement offer a second, geometry-dependent recovery route. Shared
subject effects can also produce an irreducible error floor, which can
be lowered only by acquiring genuinely additional information, as in
the ordinal-anchor construction.

The five focused simulations illustrate these phenomena under specified
synthetic designs; they do not establish comprehensive superiority, and
theory highlights potential failures of weight coverage, cutoff
calibration, and overly optimistic independence assumptions. The current results do not establish general
superiority over fitted ranking mixtures, recovery under arbitrary
partial rankings or ties, or validity of learned neural and bootstrap
candidate generators. Future work requires data-adaptive yet
independently validated calibration, quantitative candidate-coverage
guarantees for learned weights, and broader empirical comparisons.

\section*{Acknowledgments}
AI tools assisted with manuscript drafting, mathematical checks, and
simulation code. The authors are responsible for the final content.

\printbibliography
\end{document}